\documentstyle[12pt]{article}

\newtheorem{Proposition}{Proposition}
\begin{document}

\begin{center}
{ \large \bf 
Representations of the deformed 
U(su(2)) and U(osp(1,2))
algebras
\footnote{Talk given by
C. Daskaloyannis, at the Barut memorial conference,
Edirne~-~Turkey, Dec. 1995}
}

\bigskip

 Dennis BONATSOS\\
{\em ECT$^*$, Villa Tambosi, Strada delle Tabarelle 286\\
 I-38050 Villazzano (Trento), Italy}\\
C. DASKALOYANNIS\\
{\em Department of Physics, Aristotle University of Thessaloniki\\
GR-54006 Thessaloniki, Greece}\\
P. KOLOKOTRONIS, and D. LENIS\\
{\em Institute of Nuclear Physics, NCSR ``Demokritos'' \\
GR-15310 Aghia Paraskevi, Attiki, Greece}\\

\vfil
{\bf Abstract}\\
\bigskip
\begin{minipage}{5.0in}{\sf
The polynomial deformations of the Witten extensions 
of the  U(su(2)) and U(osp(1,2)) algebras are three generator algebras
with normal ordering, admitting a two generator subalgebra.  
The modules and the representations of these algebras 
are based on the construction of Verma modules,
 which are quotient modules, generated by  ideals of the original 
algebra. This construction unifies a large number of the known algebras
under the same scheme. The finite dimensional representations show new 
features such as the multiplicity of representations 
of the same dimensionality,
or the existence of finite dimensional representations only for 
some dimensions.}
\end{minipage}
\end{center}

\vfil

\section{Introduction}
Quantum algebras \cite{Drinfeld} are recently attracting much attention,
both because of their rich mathematical structure and their potential
applications in physics. The introduction of the q-deformed harmonic
oscillator \cite{Biedenharn,Macfarlane,SunFu} has been soon thereafter
followed by the introduction of the generalized deformed oscillator 
\cite{Dask1,SUSY,ArikCelik}. The idea of the generalized oscillator 
was also
introduced in mathematics twenty  years ago,  
under a totally different context  and
using a different language by  Joseph\cite{Joseph} one year
after the introduction of the Q- oscillator algebra by Arik 
and Coon\cite{ArikCoon} 

In a similar way the study of the quantum algebras su$_q$(2) 
and su$_q$(1,1) \cite{Biedenharn,Macfarlane,KD} has been followed by the
introduction of generalized deformed versions of su(2) and su(1,1). The
first generalized versions of su(2) in physics were introduced by Ro\v{c}ek 
\cite{Rocek} and Polychronakos \cite{Polychronakos}.
It should be noticed
that the same problem has been independently studied (using a slightly
different language) in a mathematical framework by Smith \cite{Smith}.
The
``generalized versions of su(2)'' of physics publications 
\cite{Rocek,Polychronakos} are
termed in mathematical publications \cite{Smith} 
as ``algebras similar to the enveloping algebra
of sl(2)''. The motivation of the physics investigations was initially 
the need to
find new models for the description of the symmetries of various physical
systems, while in mathematics the goal was mainly the construction of new
examples of noncommutative noetherian rings (see also related references in
Lenczewski\cite{Lenczewski}).

A more general framework for nonlinear deformations of U(su(2)) and 
U(su(1,1))
 has been introduced by Delbecq and Quesne \cite{DQ1,DQ2,DQ3} as 
a  generalization of the Witten\cite{Witten} algebras.  

The study of the representation theory of these nonlinear algebras is in
general an open problem, presenting features which  are not present in the
case of the corresponding Lie algebras or superalgebras. The main new
features are the possibility for the existence of multiple modules of the
same dimensionality, as well as the selective existence of modules for some
dimensionalities. These variations reveal additional difficulties on the
realization of the Hopf algebra structure at these algebras.

Sun and Li\cite{Sun} studied
 the representation theory of the deformed U(su(2)) algebra of
Ro\v{c}ek\cite{Rocek} and Polychronakos\cite{Polychronakos},
 while from the rich
variety of deformations introduced by Delbecq and Quesne\cite{DQ1}
 only the special cases
of the algebras ${\cal A}(2,1)$ \cite{DQ2} and ${\cal A}(3,1)$ \cite{DQ3}
have been considered in some detail.  In parallel, special deformations of
U(su(2)) and U(su(1,1))
 possessing representation theories as close as possible to
the usual U(su(2)) and U(su(1,1))
 algebras have been introduced and their possible
applications in physics considered \cite{Kolokotronis,Pan}.

The  study of the representations of the polynomial
deformations of the U(sl(2)) (or U(su(2)) or U(su(1,1))),
 given by Smith\cite{Smith},
is based on the construction of Verma modules which are quotient modules,
generated by simple ideals of the original algebra. This construction
was applied independently by Sun and Li\cite{Sun} at the same algebras.
This method was applied to the polynomial deformations of
U(osp(1,2)) by Van der Jeugt and Jagannathan\cite{JJ}. 
The same method 
(but using different terminology) was previously 
applied by Gruber {\em et al.}\cite{GK} to the U(su(2)) algebra and by 
Doebner {\em et al.}\cite{DGL} to the U(su(1,1)) algebra, and by
Gruber and Smirnov\cite{GS}
 to the deformed  U${}_q$(su(2)) and U${}_q$(su(1,1)).
The scope of the present presentation
 is to give a general method for the study of
the representation theory of all nonlinear deformed U(su(2)) and U(su(1,1))
algebras introduced by Delbecq and Quesne\cite{DQ1,DQ2,DQ3},
 and to their U(osp(1,2)) 
counterparts  by appropriately generalizing the method
used by Smith\cite{Smith} and Van der Jeugt and Jagannathan\cite{JJ}.

\section{ Verma modules of the deformed su(2) algebras}

Let us consider the enveloping algebra ${\cal A}$ 
\begin{equation}
{\cal A}= \mbox{\bf C}\big[ J_+, J_-, J_0],  \label{eq:defA}
\end{equation}
generated by all the polynomial combinations of the generators $J_\pm, J_0, 1
$, which satisfy the relations: 
\begin{eqnarray}
J_0 J_+ = J_+ G(J_0),  \label{eq:alg1} \\
J_- J_0 = G(J_0) J_- ,  \label{eq:alg2} \\
J_- J_+ = s J_+ J_- + f(J_0 ),  \label{eq:alg3}
\end{eqnarray}
where $s$ is a complex number and $G(z)$, $f(z)$ are polynomials of order $%
\lambda$ and $\mu$ respectively: 
\begin{equation}
G(z) = \sum\limits_{n=0}^{\lambda} g_n z^n, \quad f(z) =
\sum\limits_{n=0}^{\mu} f_n z^n.  \label{eq:sf}
\end{equation}
In the special case of $s=1$ the algebras are deformed generalizations of
the U(sl(2))~\cite{Smith} (~or U(su(2)) or U(su(1,1))~\cite
{Rocek,Polychronakos}) algebras. In the case of $s=-1$ the above defined
algebras are generalizations of the U(osp(1,2)) algebra \cite{JJ}.

The existence of the rules (\ref{eq:alg1}), (\ref{eq:alg2}) and (\ref
{eq:alg3}) imply that the algebra {\cal A} accepts a normal ordering $J_+
\prec J_0 \prec J_-$. (In physics the notion of the normal ordering
corresponds to the notion of the lexicographic order used in the
mathematical publications.) This means that: 
\begin{Proposition}\label{P.ordering}
The algebra ${\cal A}$  is a vector space with a basis of the monomials:  
\begin{equation}
\left(J_+\right)^m \left(J_0\right)^n \left(J_-\right)^p,
\qquad m,n,p \in {\bf N}, 
\label{eq:base1}
\end{equation}
or
$$
{\cal A} = 
{\rm Span}\left\{ \left(J_+\right)^m \left(J_0\right)^n \left(J_-\right)^p
\, \vert\, m,\, n,\, p \in {\bf N}
\right\}, 
$$
\end{Proposition}
where ${\bf N}$ is the set of all natural numbers plus the zero.

The rule (\ref{eq:alg1}) means that the algebra: 
\[
{\cal B}_{+} = {\bf C}\left[ J_+,\, J_0\right] 
\]
is a subalgebra of the original algebra ${\cal A}$ with normal ordering,
i.e. 
\[
{\cal B}_{+} = {\rm Span}\left\{ \left(J_+\right)^m \left(J_0\right)^n \,
\vert\, m,\,n \in {\bf N} \right\}. 
\]

Let us define the left principal \cite{Bourbaki1}[ch. I, \S 8.7] ideal 
\begin{equation}
{\cal W}_{+} = {\cal A} J_- = \left( J_- \right),  \label{eq:ideal1}
\end{equation}
and the quotient space 
\begin{equation}
{\cal T}_{+} = {\cal A} / {\cal W}_{+} .  \label{eq:ideal2}
\end{equation}
The regular projection $\pi_+$ 
\[
{\cal A} \stackrel{\pi_+}{\longrightarrow }{\cal T}_{+} 
\]
is defined by: 
\[
\pi_+(A)= A \; \mbox{mod } {\cal W}_{+}, \qquad \forall A\in {\cal A} . 
\]
The base of the vector space ${\cal T}_{+}$ is given by: 
\begin{equation}
F(m,n) = \pi_+ \left( \left(J_+\right)^m \left(J_0\right)^n \right) \, 
\stackrel{\mbox{def}}{=}\, \left(J_+\right)^m \left(J_0\right)^n \mbox{ mod }
{\cal W}_{+} .  \label{eq:base2}
\end{equation}
There is a bijection from the vector space ${\cal T}_{+}$ onto the
subalgebra ${\cal B}_{+}\subset {\cal A}$, defined as follows: 
\begin{equation}
{\cal T}_{+}\;\stackrel{i}{\longleftrightarrow}\; {\cal B}_{+}={\bf C} \big[ %
J_+,J_0\big],  \label{eq:defB}
\end{equation}
and 
\[
{\cal T}_{+} \ni \pi_+( J_+^mJ_0^n) \stackrel{i}{\longleftrightarrow }%
J_+^mJ_0^n \in {\cal B}_{+} . 
\]
Let $A\in {\cal A}$, $b=\pi_+(B)$, and $\mu_+$ be the mapping: 
\[
{\cal A}\ni A \stackrel{\mu_+}{\longrightarrow }\mu_+(A)\in 
\mbox{End}({\cal T}_{+}), 
\]
where $\mbox{End}({\cal T}_{+})$ is the set of the linear transformations
defined on ${\cal T}_{+}$. This mapping is given by: 
\[
\mu_+(A) b = \pi_+ (A \pi_+^{-1}(b))= \pi_+ (AB). 
\]
We can verify indeed that the above mapping is a left module, which is
called quotient module \cite{Bourbaki1}[ch. II, \S 1.3, example(6)] of the
algebra ${\cal A}$. This fact can be verified by calculating the action of
the generators on the base (\ref{eq:base2}). This procedure, however, leads
to complicated equations not necessary for our purposes.

The next step is to consider the left principal ideal of the algebra ${\cal B%
}_{+}$ 
\begin{equation}
{\cal W}_{+} (\eta) = {\cal B}_{+} ( J_0 - \eta ) ,  \label{eq:ideal3}
\end{equation}
and the quotient space 
\begin{equation}
{\cal T}_{+} (\eta) = {\cal B}_{+} / {\cal W}_{+}(\eta) .  \label{eq:ideal4}
\end{equation}
In a way similar to the previous one we define the regular projection $%
\pi_{\eta}$ 
\[
\begin{array}{ccc}
{\cal A} & \stackrel{\pi_+}{\longrightarrow} & {\cal T}_{+} \\ 
{\scriptstyle\pi}{\Big\downarrow} & \; & {\Big\updownarrow} {\scriptstyle i}
\\ 
{\cal T}_{+} (\eta) & \stackrel{\pi_{\eta}}{\longleftarrow} & {\cal B}_{+}
\end{array}
\]
Let us define: 
\[
\pi=\pi _{\eta} \circ {i} \circ \pi_+ . 
\]
The base of the vector space ${\cal T}_{+}(\eta)$ is given by: 
\begin{equation}
F(\eta,m) =\pi \left( \left(J_+\right)^m \right) \, \stackrel{\mbox{def}}{=}%
\, \left(J_+\right)^m \mbox{ mod } {\cal W}_{+}(\eta) .  \label{eq:base3}
\end{equation}
The vector space ${\cal T}_{+}(\eta)$, which is a quotient module of the
algebra ${\cal B}_{+}$, defines a left module of the original enveloping
algebra ${\cal A}$. Let $A\in {\cal A}$, $b=\pi(B)$, and $\mu$ be the
mapping: 
\[
{\cal A}\ni A \stackrel{\mu}{\longrightarrow }\mu(A)\in \mbox{End}({\cal T}%
_{+}(\eta)), 
\]
which is given by: 
\[
\mu(A) b = \pi (A \pi^{-1}(b))= \pi (AB) . 
\]
The above mapping defines a left module of the algebra ${\cal A}$. We must
point out that the vector space ${\cal T}_{+}(\eta)$ is isomorphic to the
subalgebra ${\cal C}\subset {\cal B}_{+} \subset {\cal A}$ defined by the
base $J_+^m$: 
\[
{\cal T}_{+}(\eta) \; \stackrel{j}{\longleftrightarrow }\; {\cal C}= %
{\bf C}\left[J_+\right]. 
\]
After calculating the action of the generators on the base (\ref{eq:base3})
we can prove the following proposition:

\begin{Proposition}\label{P.module1}
The enveloping algebra {\cal A}, the generators of which  satisfy
the following relations:
\begin{eqnarray*}
J_0 J_+ = J_+ G(J_0) , \\
J_- J_0  =  G(J_0) J_- , \\
J_- J_+  = s J_+ J_- + f(J_0 ) , 
\end{eqnarray*}
has a left module ${\cal T}_{+}(\eta)$ such that:
\begin{equation}
\begin{array}{l}
\mu(J_+) F(\eta,m) = F( \eta,m+1) , \\
\mu(J_0) F(\eta,m) = G^{[m]}(\eta) F( \eta,m) , \\
\mu(J_-) F(\eta,m) = \Phi ( \eta, m) F( \eta, m-1) , 
\end{array}
\label{eq:module1}
\end{equation}
\end{Proposition}
where: 
\begin{equation}
G^{[m]} = \underbrace{G\circ G \circ G \circ \ldots \circ G}_{m \mbox(%
times)}, \quad G^{[0]}(z) = \mbox{Id}(z)=z ,  \label{eq:defG}
\end{equation}
and 
\begin{equation}
\Phi ( \eta, 0)=0, \quad \Phi(\eta,m)= \sum\limits_{k=1}^{m} s^{k-1} f\left(
G^{[m-k]}(\eta) \right) .  \label{eq:defPhi}
\end{equation}

This representation corresponds to a ``minimum weight'' module, because: 
\begin{equation}
\mu(J_-) F(\eta,0)=0.  \label{eq:mimimum}
\end{equation}
Using the same methods as in ref. \cite{DQ1} the ``Casimir'' operator can be
defined: 
\begin{equation}
C= J_+J_- + \rho(J_0)
\label{eq:Casimir1} 
\end{equation}
where the function $\rho(z)$ satisfies the ``consistency'' equation: 
\begin{equation}
s \rho(z)- \rho(G(z))= f(z).  \label{eq:cons1}
\end{equation}
The ``Casimir'' operator satisfies the following equations: 
\[
\big[ C, J_0\big]=0, \qquad \big[J_- , C\big]_{s}=0, \qquad \big[C , J_+\big]%
_{s}= 0, 
\]
where 
\[
\big[A , B]_{q}= AB - q BA . 
\]
We can show the following relation: 
\begin{equation}
\mu(C) F (\eta,m)= s^m \rho(\eta) F(\eta,m) .  \label{eq:CasimirEigenvalues}
\end{equation}
This is the representation in which both $J_0$ and the ``Casimir'' operator $%
C$ are diagonal.

In the special case of $s=1$ the algebras are deformed generalizations of
the U(sl(2)) (~or U(su(2))  or U(su(1,1))~) algebra, the operator $C$ being
indeed the Casimir operator of the algebra. In the case of $s=-1$ the above
defined algebras are generalizations of the U(osp(1,2)) algebra, the
operator $C$ satisfying the relations: 
\[
\big[ C, J_0\big]=0, \quad \big\{ C , J_- \big\}= 0, \quad \big\{ C , J_+%
\big\}= 0. 
\]
One can then define the operator 
\begin{equation}
C_2= C^2= \left( J_+ J_- + \rho (J_o ) \right)^2,  \label{eq:Casimir2}
\end{equation}
which is the Casimir operator of the algebra: 
\[
\big[ C_2, J_0\big]=0, \quad \big[  J_-, C_2\big]= 0, \quad \big[C_2 , J_+%
\big]= 0. 
\]

\begin{Proposition}\label{P.Casimir}
In the general case $s= {\rm e}^{2 i \pi/k}$, $k\in {\bf N}$,
 we can  find a Casimir operator
${\cal D}$, which is a function of the   
operator $C$:
\begin{equation}
{\cal D} =  C^k 
\label{eq:TrigonometricCasimir}
\end{equation}
The eigenvalues of the above  Casimir are 
calculated using eqs (\ref{eq:CasimirEigenvalues}) and
(\ref{eq:TrigonometricCasimir})
$$
\mu\left({\cal D}\right) F\left(\eta , m \right)=
 \rho^k(\eta)   F\left(\eta , m \right).
$$
\end{Proposition}
Using the properties of the dual algebra ${\cal A}^\ast$ we can prove
the following proposition:
\begin{Proposition}
Let $\eta_N$ be a root of the function $\Phi(\eta,N)$ for some
natural number $N$ such that:
\begin{equation}
\Phi(\eta_N,N)\,=\, 0,
\label{eq:BasicEquation}
\end{equation}
and
$$
\Phi(\eta_N,m)\,  \ne \, 0,
\quad \mbox{for}\quad   m=1,2,\ldots, N-1
$$
then the ${\cal T}_{+}(\eta_N)$ is a N-dimensional
left   ${\cal A}$-module and its dual 
${\cal T}^\ast(\eta_N)$ is a right ${\cal A}$-module.
In this case the  $\mu(J_0)$ defined by equation (\ref{eq:module1})
is a $N\times N$ matrix satisfying the Cayley-Hamilton equation:
\begin{equation}
\mu\left(
{\left(J_0 -\eta\right) \cdot
\left( J_0 - G(\eta_N) \right)
\cdot
\left( J_0 - G^{[2]}(\eta_N) \right)
\cdots
\left( J_0 - G^{[N-1]}(\eta_N) \right)}\right) \, = \, 0 
\label{eq:CH}
\end{equation} 
and
\begin{equation}
\mu\left( J_\pm^N \right) =0
\label{eq:CH1}
\end{equation}
\end{Proposition} 

Each algebra is then characterized by the functions $G(z)$, $f(z)$ and the
constant $s$. Starting from these elements we can construct the
representation of the algebra and we can calculate the functions
$\Phi(\eta, m)$ and $\rho(\eta)$ and we find a normalized basis:
\begin{equation}
\phi (\eta ,m)=\frac 1{\sqrt{[\eta ,m]!}}F(\eta ,m),  \label{eq:normbase}
\end{equation}
where: 
$$
\begin{array}{c}
\left[\eta,m\right]= \left| \Phi (\eta ,m)\right|=
\left|\, \left[\left[ \eta, m \right]\right]\, \right|\\
\left[\eta ,0\right]!=1,
\quad \left[\eta ,m\right]!=\left[\eta,m\right] 
\left[\eta ,m-1\right]!.
\end{array}
$$
The ``normalized'' basis satisfies the following relations: 
\begin{equation}
\begin{array}{l}
\mu (J_{+})\phi (\eta ,m)=
\sqrt{\left[ \eta ,m+1\right] }\phi (\eta
,m+1), \\ 
\mu (J_0)\phi (\eta ,m)=G^{[m]}(\eta )\phi (\eta ,m),\quad \mu (C)\phi (\eta
,m)=s^m\rho (\eta )\phi (\eta ,m), \\ 
\mu (J_{-})\phi (\eta ,m)={\rm sign}\left( \Phi (\eta ,m)\right) 
\sqrt{\left[\eta ,m\right] }\phi (\eta, m-1 )=\\
=\displaystyle
\frac{\left|\, \left[\left[ \eta, m \right]\right]\, \right|}
{\sqrt{\left[ \eta ,m\right] }}\phi (\eta, m-1 )
,
\end{array}
\label{eq:module2}
\end{equation}
where ${\rm sign}\left( x\right) $ is the sign of the number $x$. 
In the  Table 1, we give examples of the functions charactering several
known algebras:

\begin{table}[thb]
\caption{Characteristic functions of the several algebras.}

\small
{
\begin{tabular}{|l|c|c|c|c|}
\hline
algebra & $G(z)$ & $s$ & $f(z)$ & $\Phi(\eta,m)$ \\ \hline\hline
U$_{q}$(su(2)) & $z+1$ & 1 & $-[2z]$ & $\left[ m\right ] \left[ -2 \eta - m
+1 \right]$ \\ \hline
U${}_q$(su(1,1)) & $z+1$ & $1$ & $[2z]$ & $\left[ m\right ] \left[ 2 \eta +
m -1 \right]$ \\ \hline
U${}_q$(osp(1 $\vert $ 2)) & $z+1/2$ & $-1$ & $-\frac{1}{4}[2z]$ & $
\begin{array}{l}
\frac{q^{\frac{1}{2}}}{4(1+q)} \Big( (-1)^{m} \left[ q^{2\eta - \frac{1}{2}
} \right] \\ 
-\left[ q^{2\eta +{m}- \frac{1}{2} } \right] \Big)
\end{array}
$ \\ \hline
${\cal A}(2, 1)$ \cite{DQ2} & $q z -1$ & $1$ & $2z (1+(1-q)z)$ & $
\begin{array}{l}
\left( 1-q^{2m} \right) \left( \eta + \frac{1}{1-q} \right) \\ 
\left( \eta - \frac{ q + q^2 + \cdots + q^{m-1} }{ 1+ q^m} \right)
\end{array}
$ \\ \hline
${\cal A}^{+}(3,1)$ \cite{DQ3} & $qz-1$ & $1$ & $2z ( 1-(1-q)^2 z^2)$ & $
\begin{array}{l}
\rho(\eta)- \\ 
\rho\left( q^m\eta - \frac{ 1 - q^m}{1-q} \right)
\end{array}
$ \\ \hline
\begin{tabular}{l}
deformed \\ 
U(osp(1, $\vert $ 2))
\end{tabular}
& $1+z $ & $-1$ & $f(z)$=polynomial & $\Phi(\eta,m)$ \\ \hline
W$^{(2)}_3$\cite{Tjin1,Tjin2} & $2+z$ & $1$ & $-(z^2+c)$ & $
\begin{array}{l}
- m\eta^2 - 2m ( m-1) \eta - \\ 
m\Big( \frac{4}{3}m^2 -2m +c + \frac{2}{3} \Big)
\end{array}
$ \\ \hline
\begin{tabular}{l}
deformed \\ 
U(su(2))\cite{Kolokotronis}
\end{tabular}
& $1+z$ & $1$ & $\phi( z(z-1) ) - \phi(z (z+1))$ & $
\begin{array}{l}
\phi((\eta+m)(\eta+m-1)) - \\ 
\phi(\eta(\eta-1) )
\end{array}
$ \\ \hline
\begin{tabular}{l}
polynomial \\ 
sl(2) \cite{Smith}
\end{tabular}
& $1+z$ & $1$ & $z^n - (z+1)^n$ & $(\eta +m )^n - \eta^n$ \\ \hline
\end{tabular}
}
\end{table}  

\section{Conclusions}\label{section:conclusions}

In conclusion, we have demonstrated that it is possible to construct a
unified representation theory of the nonlinear deformations of the su(2),
su(1,1), osp( 1$\vert$ 2) and sl(2) algebras. The proposed construction
permits the calculation of the representations of a large number of 
three generator algebras, admitting a normal order and having a 
two generator subalgebra. The deformed algebras show 
features which can not be seen
in the ordinary or q-deformed algebras. We have shown that these 
algebras have Casimir like operators, obeying  deformed commutation
relations with the generators of the algebra. The main result of this
study is that the construction of the Verma modules, and hence the
construction of the respective representations, is facilitated
 by the existence of a chain of subalgebras with normal ordering.
 The notion of the normal ordering
could be the leading idea for generalisations of this construction
in the cases of algebras with more than three generators.


\begin{thebibliography}{99}

\bibitem{Drinfeld}  V. G. Drinfeld, in {\it Proceedings of the International
Congress of Mathematicians}, ed. A. M. Gleason (American Mathematical
Society, Providence, RI, 1986) p. 798.

\bibitem{Biedenharn}  L. C. Biedenharn, J. Phys. A 22 (1989) L873.

\bibitem{Macfarlane}  A. J. Macfarlane, J. Phys. A 22 (1989) 4581.

\bibitem{SunFu}  C. P. Sun and H. C. Fu, J. Phys. A 22 (1989) L983.

\bibitem{Dask1}  C. Daskaloyannis, J. Phys. A 24 (1991) L789.

\bibitem{SUSY}  D. Bonatsos and C. Daskaloyannis, Phys. Lett. B 307 (1993)
100.

\bibitem{ArikCelik} 
M. Arik and S. Celik, Z. Phys. C59 (1993) 99.
\bibitem{Joseph} A. Joseph, Isr. J. of. Math. { 28} (1977) 177.

\bibitem{ArikCoon} M. Arik and D.D. Coon, J. Math. Phys. {17} (1976) 524.

\bibitem{KD}  P. P. Kulish and E. V. Damaskinsky, J. Phys. A 23 (1990) L415.

\bibitem{Rocek}  M. Ro\v{c}ek, Phys. Lett. B 255 (1991) 554.

\bibitem{Polychronakos}  A. P. Polychronakos, Mod. Phys. Lett. A 5 (1990)
2325.

\bibitem{Smith}  S.~P.~Smith, Trans. Am. Math. Soc. 322 (1990) 285.

\bibitem{Lenczewski}
R.~Lenczewski,  in {\em Symmetries in Science V}, ed.  B. Gruber 
{\em et al.}, (Plenum Press, NY, 1991) p. 447.


\bibitem{DQ1}  C.~Delbecq and C.~Quesne, J. Phys. A 26 (1993) L127.

\bibitem{DQ2}  C. Delbecq and C. Quesne, Phys. Lett. B 300 (1993) 227.

\bibitem{DQ3}  C. Delbecq and C. Quesne, Mod. Phys. Lett. A 8 (1993) 961.

\bibitem{Witten} E. Witten, Nucl. Phys. B330 (1990) 285.

\bibitem{Sun}  C. P. Sun and W. Li, Commun. Theor. Phys. 19 (1993) 191.

\bibitem{Kolokotronis}  D. Bonatsos, C. Daskaloyannis and P. Kolokotronis,
J. Phys. A 26 (1993) L871.

\bibitem{Pan}  F. Pan, J. Math. Phys. { 34} (1994) 5065.

\bibitem{JJ}  J. Van der Jeugt and R. Jagannathan, J. Math. Phys. 36 (1995)
4507.


\bibitem{GK}  B. Gruber and A. Klimyk, J. Math. Phys. {25} (1984) 755; 
B. Gruber, R. Lenczewski and M. Lorente, J. Math. Phys. {31} (1990) 587;

\bibitem{DGL} H. D. Doebner, B. Gruber and M. Lorente, J. Math. Phys.
{30} (1989) 587

\bibitem{GS} B. Gruber and Yu. F. Smirnov, in 
{\em Symmetries in Science V}, ed.  B. Gruber {\em et al.}, 
 (Plenum Press, NY, 1991), p. 293
 

\bibitem{Bourbaki1}  N. Bourbaki, {\it Algebra I} (Addison--Wesley,
Reading, 1974).



\bibitem{Tjin1}  T. Tjin, Phys. Lett. B 292 (1992) 60.

\bibitem{Tjin2}  J. de Boer and T. Tjin, Commun. Math. Phys. 158 (1993) 485.



\end{thebibliography}
\end{document}